\documentclass[12pt]{article}
\usepackage[utf8]{inputenc}
\usepackage[T1]{fontenc}
\usepackage{amsmath}
\usepackage{amsfonts}
\usepackage{amssymb}
\usepackage[english]{babel}
\usepackage{tikz-cd}
\title{A proof of Cairns' conjecture
about geodesic triangulations of the 2-sphere}
\date{}

\author{Jean Cerf}

\begin{document}

\maketitle
\vspace{-2.5 cm}
\begin{center}
 Dedicated to Jean-Pierre Serre
{\section*{Introduction}}
\end{center}
In a short abstract (\cite{cairns1}, 1941) S. Cairns announced the following result: «\,The space of geodesic triangulations of $S^2$ with an $n$-vertices simplicial sphere as model is homeomorphic to $ \\ SO(3)\times \mathbb{R}^{2n-3}$\,». In his 1943 Annals paper (cf \cite{cairns2}) Cairns proved the connexity of that space, a result we shall call «\,Cairns' theorem\,», the initial statement being called «\,Cairns' conjecture\,».
Cairns' theorem is one of the pillars of the present paper as well as two classical planar results:
\begin{enumerate}
    \item the BCH theorem (Bloch-Connelly-Henderson \cite{bloch}, 1981):\newline «\,The space of SL triangulations of a convex SL disk in $\mathbb{R}^2$ modeled on a simplicial disk with $n$ interior vertices is homeomorphic to $\mathbb{R}^{2n}$\,».
    \item the theorem of Wagner (1936) and Fáry (1948) about linearisation of planar graphs (cf \cite{richter}).
\end{enumerate}
The proof itself is based on the study of the space $\mathcal{K}$ as a differentiable manifold, whose dimension grows linearly with the number of vertices of $K$. The notions and tools we use belong to differential topology, including Lie groups actions, stratifications, fibrations. A possible similar rewriting of BCH theorem's proof would lead to a stronger final statement: the space $\mathcal{K}$ is diffeomorphic to $SO(3) \times \mathbb{R}^{2n-3}$.

Thanks for their support to my friends François Laudenbach, Alexis Marin and Valentin Poenaru.

\subsection*{ 1 Cairns' conjecture.}

The sphere $S^2$ being defined as the sphere with center $O$, radius 1, in euclidean
$\mathbb{R}^3$, a geodesic 2-simplex of $S^2$ is the radial projection from $O$ of an affine triangle
in $\mathbb{R}^3$ whose plane does not contain $O$. A geodesic 1-simplex is an arc of a great circle
of length $< \pi$. A triangulation of $S^2$ is said to be \textbf{simplicially geodesic} if all its simplices
are geodesic. The sphere $S^2$ equipped with such a triangulation is called an \textbf{SG sphere}.

The SG spheres are divided into families according to their respective simplicial models. We agree to
take as a model for each family a positively oriented SG sphere (i.e., all 2-simplices are oriented by
the positive orientation of $S^2$). We denote this sphere by $K$ and we number its vertices from 1
to $n$.

A positively oriented SG 2-sphere, denoted $K'$, is said to be \textbf{modeled by $K$} if one has
chosen a bijection of the 0-skeletons: $K^0 \to K'^0$ compatible with the simplicial structures and
the orientations. The numbering of the vertices of $K$ determines that of $K'$, which allows
identifying $K'$ with a point of $(S^2)^n$. The union of these points for all $K'$ is an open set of
$(S^2)^n$ whose differential structure, and thus the topology, are independent of the numbering
chosen for the vertices of $K$. Equipped with this topology, this open set is called the \textbf{space
of SG spheres modeled by $K$} and denoted by $\mathcal{K}$.

\begin{center}
\textbf{Cairns' Conjecture.} There exists a homeomorphism
\begin{equation*}
\mathcal{K} \approx SO(3) \times \mathbb{R}^{2n-3}
\end{equation*}
\end{center}

This formula is easily verified for the boundary of the tetrahedron ($n=4$).

\subsection*{ 2 \textbf{The stratified space of SG triangles of $S^2$}.}

An \textbf{SG triangle} is an SG circle with three 1-simplices.
The space of ordered SG triangles (i.e. the vertices are numbered 1, 2, 3) is identified with an open
set of $(S^2)^3$; it will be convenient to use, besides its topology, its differentiable structure; it will
be denoted by $\mathcal{B}$.

The space $\mathcal{B}$ is the disjoint union
\begin{equation*}
\mathcal{B} = \mathcal{B}_+ \cup \mathcal{B}_- \cup \mathcal{B}_0
\end{equation*}
where the elements of $\mathcal{B}_+$ and $\mathcal{B}_-$ are the boundaries of ordered 2-
simplices, classified as\\ + or - according to their orientation. The elements of $\mathcal{B}_0$ are the
flat ordered triangles (i.e. these contained in a great circle) which have no privileged orientation.


 One denotes by $DSL(3)$ the group of diagonal $3 \times 3$ matrices with positive entries and determinant 1.
As is shown by its representation in $\mathbb{R}^3$, $DSL(3)$ is homeomorphic to $\mathbb{R}^2$.

The group $SL(3)$ acts on $S^2$ by composition of the action on the ambient $\mathbb{R}^3$, followed by the radial projection of $\mathbb{R}^3 \setminus \{0\}$ on $S^2$. That action extends to $(S^2)^3$.

\textbf{Lemma} \begin{enumerate}\item \emph{The group $SL(3)$ acts transitively on $\mathcal{B}_+$, $\mathcal{B}_-$ and $\mathcal{B}_0$.}

\item \emph{The restriction of the action of $SL(3)$ on any triple $(b_+, b_-, b_0) \in (\mathcal{B}_+, \mathcal{B}_-, \mathcal{B}_0)$ admits continuous sections over $\mathcal{B}_+$, $\mathcal{B}_-$, $\mathcal{B}_0$ and leads therefore to product decompositions of $SL(3)$ over $\mathcal{B}_+$, $\mathcal{B}_-$, $\mathcal{B}_0$ respectively.}

\item \emph{In the case of $\mathcal{B}_+$ and $\mathcal{B}_-$, the fibers of the preceding decomposition are homeomorphic to $DSL(3)$. There exists for instance an homeomorphism $SL(3) \approx \mathcal{B}_+ \times DSL(3)$.}
\end{enumerate}
The proof of the lemma reduces to 2. and 3.

\underline{Proof of 2.} The case of $\mathcal{B}_+$ and $\mathcal{B}_-$ being similar, we restrict ourselves to $\mathcal{B}_+$ and $\mathcal{B}_0$.

\underline{Case of $\mathcal{B}_+$.}  Choosing any $b_+ \in \mathcal{B}_+$, we have only to exhibit a section for the $SL(3) \to \mathcal{B}_+$ map: $g \mapsto g \cdot b_+$. We take for instance the map
    \[
    b \longmapsto \frac{g(b)}{\det g(b)}
    \]
    where $g(b)$ denotes the element of $GL_+(3)$ such that $(g(b)) \cdot b = b_+$.

\underline{Case of $\mathcal{B}_0$}. For any $b_0 \in \mathcal{B}_0$, there exists $g \in SO(3)$ such as $g \cdot b_0$ is equatorial with positive orientation (in the sense of the horizontal $\mathbb{R}^2$). The fact that two such SG triangles are \newline $SL(3)$-equivalent is a result in elementary planar geometry.

\underline{Proof of 3.} Let $b_+$ be the point of $\mathcal{B}_+$ defined by the canonical basis of $\mathbb{R}^3$. The subgroup of $SL(3)$ fixing $b_0$ is $DSL(3)$ \hfill $\square$

\bigskip

\noindent \textbf{Volume of an SG triangle}. For all $\tau = (s_i, s_j, s_k) \in \mathcal{B}$, we set
\begin{equation*}
V(\tau) = \det(s_i, s_j, s_k),
\end{equation*}
where $s_i, s_j, s_k$ are considered as elements of $\mathbb{R}^3$. We will say that $V(\tau)$ is
the \textbf{volume} of $\tau$ (as opposed to the Riemannian area on $S^2$, which is inadequate
here).

\vspace{0.5cm}

\noindent \textbf{Expression of the volume in spherical coordinates} $(\theta, \varphi)$. \\
($\theta \in S^1 : \text{longitude}, \varphi \in \mathbb{R} : \text{latitude}$).
\begin{itemize}
\item If no vertex is at a pole:
\begin{equation*}
V(\tau) = \cos\varphi_i \cos\varphi_j \cos\varphi_k \left[ \tan\varphi_i \sin(\theta_k -
\theta_j) + \tan\varphi_j \sin(\theta_i - \theta_k) + \tan\varphi_k \sin(\theta_j - \theta_i) \right]
\end{equation*}
\item If for example $s_i$ is at the north pole:
\begin{equation*}
V(\tau) = \cos\varphi_j \cos\varphi_k \sin(\theta_k - \theta_j).
\end{equation*}
\end{itemize}

\vspace{0.5cm}

\noindent \textbf{Two properties of the volume function}.
\begin{enumerate}
\item \emph{We have: $(\mathcal{B}_+, \mathcal{B}_-, \mathcal{B}_0) = V^{-1}(\mathbb{R}_+,
\mathbb{R}_-, 0)$.}
\item \emph{$\mathcal{B}_0$ is a two-sided differentiable submanifold of codimension 1 of
$\mathcal{B}$, separating $\mathcal{B}_+$ and $\mathcal{B}_-$.With other words: $(\mathcal{B}_+, \mathcal{B}_-, \mathcal{B}_0)$ is a codimension one $DIFF$ stratification of  $\mathcal{B}$. }

\end{enumerate}

\noindent \textbf{Proof of property 2}. We prove that the function $V$ is a non-singular Morse
function on $\mathcal{B}_0 = V^{-1}(0)$. The stratified space $\mathcal{B}$ being invariant under
rotation of $S^2$, it suffices to test this non-singularity for a triangle $\tau$ lying on the equator.
Now, for such a $\tau$, we have for instance

\begin{equation*}
\frac{\partial V}{\partial \varphi_i} = \sin(\theta_k - \theta_j) \neq 0. \qquad \square
\end{equation*}

\subsection*{ 3 \textbf{Disks with triangular  boundary: stratified fibration theorem}.}

Let $K$ and $\mathcal{K}$ be as in $\S 1$
and let $\tau$ be a 2-simplex of $K$. We denote by $D$ the
SG disk complementary to $\tau$ and by $\mathcal{D}$ the space of SG disks positively modeled by
$D$.
For any $D' \in \mathcal{D}$, the boundary $\partial D'$ modeled by $\partial D$ is an SG triangle;
we denote by $\delta : \mathcal{D} \to \mathcal{B}$ the map thus defined. The inverse image under
$\delta$ of the stratification $(\mathcal{B}_+, \mathcal{B}_-, \mathcal{B}_0)$ of $\mathcal{B}$ is
denoted $(\mathcal{D}_+, \mathcal{D}_-, \mathcal{D}_0)$. 
An element of $\mathcal{D}$ is in $\mathcal{D}_+$, $\mathcal{D}_-$ or $\mathcal{D}_0$
depending on whether its support strictly contains a closed hemisphere, is contained in an open
hemisphere, or is a closed hemisphere.("Hemisphere", without specifying North or South, means
"one of the halves of $S^2$ cut out by a great circle".)

\emph{The space  $\mathcal{D}_+$ is canonically homeomorphic to $\mathcal{K}.$
 The fact that $\mathcal{D}_-$ and $\mathcal{D}_0$ are not empty is part of the following results.}

\textbf{Theorem 1.} \textbf{The map $\delta : \mathcal{D} \to \mathcal{B}$ is surjective and defines
on each of the strata $\mathcal{D}_+, \mathcal{D}_-, \mathcal{D}_0$ a trivial fibration with
respective basis $\mathcal{B}_+, \mathcal{B}_-, \mathcal{B}_0$}.

\textbf{Corollary of theorem 1.} \textbf{With the exception of the 2-simplex, every simplicial 2-disk
with triangular boundary can be realized by an SG triangulation of the closed hemisphere.}

\textbf{Proof of the corollary}: it explicits ''$\mathcal{D}_0 \neq \emptyset$ ''. $\square$

\textbf{Proof of the surjectivity of $\delta$}. We denote by $K'$ the triangulation of $S^2$ obtained
by adding to $D$ a star subdivision of $\tau$, and by $\mathcal{K}'$ the space of SG triangulations of
$S^2$ positively modeled by $K'$.

Every disk of $S^2$ bounded by an SG triangle possesses interior points with respect to which it is
geodesically star-shaped. (If it is a 2-simplex or a hemisphere we can take its barycenter; otherwise
we take the antipode of the barycenter of the complementary disk). Therefore the restriction map
$\mathcal{K}' \to \mathcal{D}$ is surjective. But $\mathcal{K}'$ is connected according to Cairns'
theorem (cf Introduction). Therefore $\mathcal{D}$ and consequently $\delta(\mathcal{D})$ are
connected.

According to 1. of the lemma, the group $SL(3)$ acts transitively in $\mathcal{B}_+$,
$\mathcal{B}_-$ and $\mathcal{B}_0$. We have therefore only to show that $\delta(\mathcal{D})$ meets
these three spaces, which is obvious for $\mathcal{B}_+$ (which contains $\delta(D)$).

To show that $\delta(\mathcal{D})$ meets $\mathcal{B}_-$, we show that $\mathcal{D}_-$ is not
empty. According to the Wagner-F\'{a}ry theorem (cf \cite{richter}, p. 19), every simplicial 2-disk can be realized
by a simplicially affine disk of $\mathbb{R}^2$. We realize $D$ by such a disk contained in a plane of
$\mathbb{R}^3 \setminus \{0\}$, which we project radially onto $S^2$. The disk $D'$ thus obtained
is an element of $\mathcal{D}_-$.

Since $\mathcal{B}_+$, $\mathcal{B}_-$ and $\mathcal{B}_0$ are characterized by the positive,
negative, zero values of the volume function (\S 2. property 1), every path joining 
 $\delta(D)$ to $\delta(D')$ in $\delta(\mathcal{D})$ meets $\mathcal{B}_0$.

At this time of the proof we can complete the notations, writing $  \delta=(\delta_+, \delta_-, \delta_0)$,  the corresponding fibers being denoted by  $(\mathcal{F}_+, \mathcal{F}_-, \mathcal{F}_0)$ .

\textbf{Proof of the three fibrations}. The surjectivity of each of the maps $\delta_+$, $\delta_-$, $\delta_0$ has been proved above. The group $\mathrm{SO}(3)$ acts in each of these maps.It means for instance that the following diagram commutes:
\[
\begin{tikzcd}
\mathrm{SO}(3) \times \mathcal{D}_+ \arrow[r] \arrow[d, "\mathrm{id} \times \delta_+"'] & \mathcal{D}_+ \arrow[d,"\delta_+"] \\
\mathrm{SO}(3) \times \mathcal{B}_+ \arrow[r] & \mathcal{B}_+
\end{tikzcd}
\]
 In such a situation, the existence, at the target level, of global sections for the group action on a chosen point implies product decomposition at the domain level ( cf \cite{cerf}, lemma 2 of the introductory chapter). Now the existence of such sections is part of 2. of the lemma.\quad$\square$

\subsection*{4. Global fibration of the space $\mathcal{D}$. Proof of the conjecture}

According to the Palais--Cerf fibration theorem (cf \cite{cerf} corollary 2, p. 294, or \cite{palais} theorem C), the boundary map $\delta : \mathcal{D} \to \mathcal{B}$ is a locally trivial fibration. The space $\mathcal{B}$ being connected, and the \emph{fibers} of the restrictions $\delta_+$, $\delta_-$, $\delta_0$ being the same as in \S~3 (only the \emph{trivialisations} are different), one has
\[
\mathcal{F}_+ \approx \mathcal{F}_- \approx \mathcal{F}_0.
\]

The space of SG triangles of the open northern hemisphere is homeomorphic, by radial projection from the origin of $\mathbb{R}^3$, to the space of simplicially linear (SL) triangles of the tangent plane to $S^2$ at the north pole, with SL-subdivisions corresponding to SG subdivisions. It follows therefore from the BCH theorem (cf \cite{bloch}) that $\mathcal{F}_-$ (and consequently $\mathcal{F}_+$) are homeomorphic to $\mathbb{R}^{2(n-3)}$.

On the other hand one has after theorem 1
\[
\mathcal{K} \approx \mathcal{D}_+ \approx \mathcal{B}_+ \times \mathcal{F}_+,
\]
and after 3. of the lemma
\[
SL(3) \approx \mathcal{B}_+ \times DSL(3) \approx \mathcal{B}_+ \times \mathbb{R}^2.
\]
Finally one has
\[
\mathcal{K} \approx \mathcal{B}_+ \times \mathbb{R}^{2(n-3)} \approx SL(3) \times \mathbb{R}^{2(n-4)} \approx SO(3) \times \mathbb{R}^{2n-3},
\]
and we can state:

\textbf{Theorem 2. [Cairns' conjecture]
Let $K$ a simplicially geodesic $2$-sphere with $n$ vertices. The space of simplicially geodesic $2$-spheres modeled by $K$ is homeomorphic to $SO(3) \times \mathbb{R}^{2n-3}$.}

\renewcommand{\refname}{Bibliography}

\end{document}